\documentclass[10pt,twoside]{art_large}
\usepackage{amssymb,amsbsy,amsmath,amsfonts,amssymb,amscd}
\usepackage[latin1]{inputenc}
\usepackage[english,francais]{babel}

\newtheorem{theorem}{Theorem}[section]
\newtheorem{lemma}[theorem]{Lemma}

\newcommand{\R}{\mathbb R}

 \ComParit{S0764-4442} \PIT{FLA} \PXMA{????} \Add{?}
 \Volume{???} \Year{2002} \FirstPage{1} \LastPage{4}
\AuteurCourant{Jean-Paul Dufour and Nguyen Tien Zung}
\TitreCourant{Nondegeneracy of the Lie algebra {$\mathfrak{aff}(n)$}} \Journal
\Rubrique{Géométrie différentielle}{Differential geometry} \SousRubrique{}{}
\PresentePar{Prénom}{NOM}
\Recu{jour mois ann\'ee}{apr\`es r\'evision}{jour mois ann\'ee}
\begin{document}
\selectlanguage{english}
\title{Nondegeneracy of the Lie algebra {$\mathfrak{aff}(n)$}}
\author{%
Jean Paul Dufour~$^{\text{a}}$ \ and \  Nguyen Tien Zung~$^{\text{b}}$}
\address{%
\begin{itemize}\labelsep=2mm\leftskip=-5mm
\item[$^{\text{a}}$]
GTA, UMR 5030 CNRS, Département de Mathématiques, Université Montpellier II\\
E-mail: dufourj@math.univ-montp2.fr
\item[$^{\text{b}}$]
Laboratoire Emile Picard, UMR 5580 CNRS, UFR MIG, Université Toulouse III \\
\end{itemize}
}
\maketitle \thispagestyle{empty}
\begin{Abstract}{%
We show that  $\mathfrak{aff}(n)$, the Lie algebra of  affine transformations of
$\R^n,$ is formally and analytically nondegenerate in the sense of A. Weinstein.
This means that every analytic (resp., formal) Poisson structure vanishing at a
point with a linear part corresponding to $\mathfrak{aff}(n)$ is locally
analytically (resp., formally) linearizable. }\end{Abstract}
\selectlanguage{french}
\begin{Ftitle}{%
Non-dégénérescence de l'algèbre de Lie $\mathfrak{aff}(n)$}\end{Ftitle}
\begin{Resume}{%
Nous montrons que toute structure de Poisson analytique (resp., formelle), qui
s'annule en un point et dont la partie linéaire correspond à  l'algèbre
$\mathfrak{aff}(n)$ des transformations affines sur ${\mathbb R}^n$, est localement
analytiquement (resp., formellement) linéarisable. }\end{Resume}

\setcounter{section}{0}
\selectlanguage{english}

\section{Introduction}

Following A. Weinstein \cite{Weinstein-Poisson1983}, we will say that a Lie algebra
$\mathfrak g$ is formally (resp. analytically, smoothly) non-degenerate if  any
formal (resp. analytic, smooth) Poisson  structure $\Pi$ vanishing at a point, with
a linear part $\Pi^{(1)}$ corresponding to $\mathfrak{g}$, is formally (resp.
analytically, smoothly) linearizable. An interesting and largely open question in
Poisson geometry is to find and classify nondegenerate Lie algebras in the above
senses. Up to now, only few nondegenerate Lie algebras are known. These include: the
nontrivial 2-dimensional Lie algebra $\mathfrak{aff}(1)$ (see appendix of
\cite{Arnold-Geometrical1988} on Poisson structures and densities); semisimple Lie
algebras  \cite{Weinstein-Poisson1983,Conn-Analytic1984,Conn-Smooth1985}; a complete
list of nondegenerate Lie algebras in dimensions 3 \cite{Dufour-Linearization1990}
and 4 \cite{Molinier-Thesis1993}, direct products of semisimple algebras with
$\mathbb R$ (or $\mathbb C$) \cite{Molinier-Thesis1993}, and  direct products of $n$
copies of $\mathfrak{aff}(1)$ \cite{DuMo-Naff1995}.

Recent works by A. Wade, Ph. Monnier and the second author on  Levi decomposition of
Poisson structures \cite{MoZu-Levi2002,Wade-Levi1997,Zung-Levi2002} open a new way
of linearizing Poisson structures by first looking for a semi-linearization
associated to a Levi decomposition of their linear part. We will recall this method
in Section 2. Using it, we obtain the following result:

\begin{theorem}
\label{theorem:affn} For any natural number $n$, the Lie algebra $ \mathfrak{aff}(n,
{\mathbb K}) = \mathfrak{gl}(n, {\mathbb K}) \ltimes {\mathbb K}^n$ of affine
transformations of ${\mathbb K}^n$, where $\mathbb K = \mathbb R$ or $\mathbb C$, is
formally and analytically nondegenerate.
\end{theorem}

The above theorem provides a new entry in a very short list of known examples of
nondegenerate Lie algebras. And unlike previously known examples,
$\mathfrak{aff}(n)$ for $n \geq 2$ is neither reductive nor solvable. We remark that
the main sisters of $\mathfrak{aff}(n)$, namely the algebra $\mathfrak{e}(n) =
\mathfrak{so}(n) \ltimes {\mathbb R}^n$ of Euclidean transformations and the algebra
$\mathfrak{saff}(n) = \mathfrak{sl}(n) \ltimes {\mathbb K}^n$ of volume preserving
affine transformations, are unfortunately degenerate for $n = 2$ or $3$, as simple
polynomial non-linearizable examples show (we suspect that they are degenerate for
$n > 2$ as well).

Another interesting feature about $\mathfrak{aff}(n)$ is that it is a Frobenius Lie
algebra, in the sense that there is a dense open subset in the dual of
$\mathfrak{aff}(n)$ where the corresponding linear Poisson structure is
nondegenerate. The existence of an open coadjoint orbit probably plays a role in the
nondegeneracy of  $\mathfrak{aff}(n)$ and some other Frobenius Lie lagebras.
However, being Frobenius does not guarantee nondegeneracy: already in dimension 4
there are counter-examples, which can be seen from the list given in
\cite{Molinier-Thesis1993}.

The rest of this Note is organized as follows: In Section 2 we obtain an improved
version of the analytic semi-linearization result of \cite{Zung-Levi2002}, which
works not only for $\mathfrak{aff}(n)$ but for many other Lie algebras as well. Then
in Section 3 we give a proof of Theorem \ref{theorem:affn} based on this improved
semi-linearization and using a trick involving Casimir functions for
$\mathfrak{gl}(n)$. For simplicity of exposition, we will restrict our attention to
the analytic case. The formal case is absolutely similar, if not simpler. In Section
4 we show that the Lie algebras $\mathfrak{saff}(2)$ and $\mathfrak{e}(3)$ are
degenerate (formally, analytically and smoothly).

\section{Semi-linearization for $\mathfrak{aff}(n)$}

Denote by
$\mathfrak{g} =\mathfrak{s}\ltimes \mathfrak{r}$
a Levi decomposition for a (real or complex) Lie algebra $\mathfrak{g}$, where
$\mathfrak{s}$ is semisimple and $\mathfrak{r}$ is the solvable radical. Let $\Pi$
be an analytic Poisson structure vanishing at a point $0$ in a manifold whose linear
part at $0$ corresponds to $\mathfrak{g}$. According to the main result of
\cite{Zung-Levi2002} (called the analytic Levi decomposition theorem), there exists
a local analytic system of coordinates $(x_1,...,x_m,y_1,...,y_d)$ in a neighborhood
of $0$, where $m = \dim \mathfrak{s}$ and $d = \dim \mathfrak{r}$, such that in
these coordinates we have
\begin{equation}
\label{eqn:Levi} \{ x_i,x_j\}=\sum c_{ij}^kx_k\ ,\ \ \{x_i,y_r\}=\sum a_{ir}^sy^s
\end{equation}
where $c_{ij}^k$ are structural constants of $\mathfrak{s}$ and $a_{ir}^s$ are
constants. This gives what we  call a semi-linearization for $\Pi.$ Note that the
remaining Poisson brackets $\{y_r,y_s\}$ are nonlinear in general.


We now restrict our attention to the case where $\mathfrak{g} = \mathfrak{aff}(n)$,
$m = n^2 -1$, $d = n +1$,  $\mathfrak{s} = \mathfrak{sl}(n), \mathfrak{r} = {\mathbb
K}({\rm Id}) \ltimes {\mathbb K}^n$ where $\rm Id$ acts on ${\mathbb K}^n$ by the
identity map.  The following lemma says that we may have a semi-linearization
associated to the decomposition $ \mathfrak{aff}(n) = \mathfrak{gl}(n) \ltimes
{\mathbb K}^n$ (which is better than the Levi decomposition).

\begin{lemma}
\label{lemma:Improvedlevi} There are local analytic coordinates $x_1,\dots ,
x_{n^2-1},y_0,y_1,\dots ,y_n$ which satisfy Relations (\ref{eqn:Levi}), with the
following extra properties: $\{ y_0,y_r\}=y_r$ for $r=1,...,n$; $\{ x_i,y_0\}=0 \
\forall i$.
\end{lemma}

{\sl Proof.} We can assume that the coordinates $y_r$ are chosen so that Relations
(\ref{eqn:Levi}) are already satisifed, and $y_0$ corresponds to $\rm Id$ in
${\mathbb K}({\rm Id}) \ltimes {\mathbb K}^n$. Then the Hamiltonian vector fields
$X_{x_i}$ are linear and form a linear action of $\mathfrak{sl}(n)$. Because of
(\ref{eqn:Levi}), we have that $\{x_i,y_0\} = 0$, which implies that
$[X_{x_i},X_{y_0}] = 0$, i.e. $X_{y_0}$ is invariant under the $\mathfrak{sl}(n)$
action. Moreover we have $X_{y_0}(x_i) = 0$ (i.e. $X_{y_0}$ does not contain
components $\partial / \partial x_i$), and $X_{y_0} = \sum_1^n y_i \partial /
\partial y_i +$  nonlinear terms. Hence we can use (the parametrized equivariant
version of) Poincaré linearization theorem to linearize $X_{y_0}$ in a
$\mathfrak{sl}(n)$-invariant way. After this linearization, we have that $X_{y_0} =
\sum_1^n y_i \partial / \partial y_i$. In other words, Relations (\ref{eqn:Levi})
are still satisfied, and moreover we have $\{ y_0,y_i\}= X_{y_0}(y_i) = y_i$.
$\square$

The significance of the above lemma is that we can extend the analytic
semi-linearization of \cite{Zung-Levi2002} from $\mathfrak{sl}(n)$ to
$\mathfrak{gl}(n) = \mathfrak{sl}(n) \oplus {\mathbb K}$ . Hereafter we will
redenote $y_0$ in the above lemma by $x_{n^2}$. Then Relations (\ref{eqn:Levi}) are
still satisfied.

{\it Remark}. Lemma \ref{lemma:Improvedlevi} still holds if we replace
$\mathfrak{aff}(n)$ by any Lie algebra of the type $(\mathfrak{s} \oplus {\mathbb
K}e_0) \ltimes \mathfrak{n}$ where $\mathfrak{s}$ is semisimple and $e_0$ acts on
$\mathfrak{n}$ by the identity map (or any matrix whose corresponding linear vector
field is nonresonant and satisfies a Diophantine condition).

\section{Linearization for $ \mathfrak{aff}(n)$}

We will work in a coordinate system $x_1,...,x_{n^2},y_1,...,y_n$ provided by Lemma
\ref{lemma:Improvedlevi}. We will fix the variables $x_1,...,x_{n^2}$, and consider
them as linear functions on $\mathfrak{gl}^\ast(n)$ (they give a Poisson projection
from our $(n^2+n)$-dimensional space to  $\mathfrak{gl}^\ast(n)$). Denote by
$F_1,\dots,F_n$ the $n$ basic Casimir functions for $\mathfrak{gl}^\ast(n)$ (if we
identify $\mathfrak{gl}(n)$ with its dual via the Killing form, then $F_1,\dots,F_n$
are basic symmetric functions of the eigenvalues of $n \times n$ matrices). We will
consider $F_1(x),\dots,F_n(x)$ as functions in our $(n^2+n)$-dimensional space,
which do not depend on variables $y_i$. Denote by $X_1,\dots,X_n$ the Hamiltonian
vector fields of $F_1,\dots,F_n$.

\begin{lemma}
\label{casimir} The vector fields $X_1, \dots, X_n$ do not contain components
$\partial /\partial x_i$. They form a system of $n$ linear commuting vector fields
on ${\mathbb K}^n$ (the space of $y=(y_1,\dots ,y_n)$) with coefficients which are
polynomial in $x=(x_1,\dots ,x_{n^2}).$  The set of $x$ such that they are linearly
dependent everywhere in ${\mathbb K}^n$ is an analytic space of complex codimension
greater than 1 (when ${\mathbb  K} = {\mathbb C}$).
\end{lemma}

{\sl Proof.} The fact that the $X_i$ are $y$-linear with $x$-polynomial coefficients
follows directly from Relations (\ref{eqn:Levi}). Since $F_i$ are Casimir functions
for $\mathfrak{gl}(n)$, we have $X_i(x_k) = \{F_i,x_k\} = 0$, and $[X_i,X_j] =
X_{\{F_i,F_j\}} = 0$.

One checks that, for a given $x$, $X_1 \wedge \dots \wedge X_n = 0$ identically on
${\mathbb K}^n$ if and only if $x$ is a singular point for the map $(F_1,...,F_n)$
from $\mathfrak{gl}^\ast(n)$ to ${\mathbb K}^n$. The set of singular points of the
map $(F_1,...,F_n)$ in the complex case is of codimension greater than 1 (in fact,
it is of codimension 3). $\square$

\begin{lemma} \label{pitilde} Write the Poisson structure $\Pi$  in the form
$\Pi =\Pi^{(1)}+\tilde\Pi$, where $\Pi^{(1)}$ is the linear part and $\tilde\Pi$
denote the higher order terms. Then $\tilde\Pi$ is a Poisson structure which can be
written in the form
\begin{equation}
\label{eqn:pitilde} \tilde\Pi =\sum_{i<j}f_{ij}X_i\wedge X_j \ ,
\end{equation}
where the functions $f_{ij}$ are analytic functions which depend only on the
variables $x$, and they are Casimir functions for $\mathfrak{gl}^\ast(n)$ (if we
consider the variables $x$ as linear functions on $\mathfrak{gl}^\ast(n)$).
\end{lemma}
{\sl Proof.} We work first locally near a point $(x,y)$ where the vector fields
$X_k$ are linearly independent point-wise. As $\tilde\Pi$ is a 2-vector field in
${\mathbb K}^n=\{ y\}$ (with coefficient depending on $x$) we have a local formula
$\tilde\Pi =\sum_{i<j}f_{ij}X_i\wedge X_j$ where $f_{ij}$ are analytic functions in
variables $(x,y).$ Since $X_k$ are Hamiltonian vector fields for $\Pi$ and also for
$\Pi^{(1)}$, we have $[X_k, \tilde\Pi] = [X_k,\Pi] - [X_k,\Pi^{(1)}] = 0$ for
$k=1,\dots,n$. This leads to $X_k(f_{ij})=0 \ \forall \ k,i,j.$ Hence, because the
$X_k$ generate ${\mathbb K}^n,$ the functions $f_{ij}$ are locally independent of
$y.$ Using analytic extension, Hartog's theorem and the fact that the set of $x$
such that $X_1,\dots,X_n$ are linearly dependent point-wise everywhere in ${\mathbb
K}^n$ is of complex codimension greater than 1, we obtain that $f_{ij}$ are local
analytic functions in a neighborhood of $0$ which depend only on the variables $x$.
The fact that $\tilde\Pi$ is a Poisson structure, i.e. $[\tilde\Pi,\tilde\Pi] = 0$,
is now evident, because $X_k(f_{ij}) = 0$ and $[X_i,X_j] = 0$.

Relations $[X_{x_k},\tilde\Pi]= [X_{x_i},\Pi] - [X_{x_i},\Pi^{(1)}] = 0$ imply that
$X_{x_k}(f_{ij})= 0$, which means that $f_{ij}$ are Casimir functions for
$\mathfrak{gl}^\ast(n)$. $\square$

{\it Remark}. Lemma \ref{pitilde} is still valid in the formal case. In fact, every
homogeneous component of $\tilde\Pi$ satisfies a relation of type
(\ref{eqn:pitilde}).

\begin{lemma} \label{primitive} There exists a vector field $Y$ of the form
$Y=\sum_{i=1}^n\alpha_i X_i,$ where the analytic functions $\alpha_i$ depend only on
the variables $x$ and are Casimir functions for $\mathfrak{gl}^\ast(n)$, such that
\begin{equation}
 [Y,\Pi^{(1)}]= - \tilde\Pi\ ,\ \ [Y,\tilde\Pi ]=0.
\end{equation}
\end{lemma}
{\sl Proof.} Since the functions $f_{ij}$ of Lemma \ref{pitilde} are analytic
Casimir functions for $\mathfrak{gl}(n),$ we have $f_{ij}=\phi_{ij}(F_1,\dots ,F_n)$
where  $\phi_{ij}(z_1,\dots ,z_n)$ are analytic functions of $n$ variables. On the
other hand, since $\Pi^{(1)},\tilde\Pi$ and $\Pi = \Pi^{(1)} + \tilde\Pi$ are
Poisson structures, they are compatible, i.e. we have $[\Pi^{(1)},\tilde\Pi] = 0$.
Decomposing this relation, we get
$ \frac{\partial\phi_{ij}}{\partial z_k}+\frac{\partial\phi_{jk}}{\partial z_i}+
\frac{\partial\phi_{ki}}{\partial z_j}=0 \ \forall \ i,j,k .$
This is equivalent to the fact that the 2-form $\phi:=\sum_{ij}\phi_{ij}dz_i\wedge
dz_j$ is closed. By Poincaré's lemma we get $\phi=d\alpha$ with an 1-form
$\alpha=\sum_i\alpha_idz_i.$ Then we put $Y:=\sum_i\alpha_i(F_1,\dots ,F_n)X_i.$ An
elementary calculation proves that $Y$ is the desired vector field. $\square$

{\it Proof of Theorem \ref{theorem:affn}}. Consider a path of Poisson structures
given by $\Pi_t:=\Pi^{(1)}+t\tilde\Pi .$ As we have $[Y,\Pi_t ]=\tilde\Pi =
\frac{d}{dt} \Pi_t$, the time-1 map of the vector field $Y$ moves $\Pi^{(1)} =
\Pi_0$ into $\Pi = \Pi_1$. This shows that $\Pi$ is locally analytically
linearizable, thus proving our theorem. $\square$

\section{Degeneracy of $\mathfrak{saff}(2)$ and $\mathfrak{e}(3)$}

The linear Poisson structure corresponding to $\mathfrak{saff}(2)$ has the form
$\Pi^{(1)}=2e\partial h\wedge\partial e-2f\partial h\wedge\partial f +h\partial
e\wedge\partial f+y_1\partial h\wedge\partial y_1-y_2\partial h\wedge\partial
y_2+y_1\partial e\wedge\partial y_2+y_2\partial f\wedge\partial y_1$ in a natural
system of coordinates. Now put $\Pi=\Pi^{(1)}+\tilde\Pi$ with
$\tilde\Pi=(h^2+4ef)\partial y_1\wedge\partial y_2.$ Then $\Pi$ is a Poisson
structure, vanishing at the origin, with a linear part corresponding to
$\mathfrak{saff}(2).$ For $\Pi^{(1)}$ the set where the rank is less or equal to 2
is a codimension 2 subspace (given by the equations $y_1=0$ and $y_2=0$). For $\Pi$
the set where the rank is less or equal to 2 is a 2-dimensional cone (the cone given
by the equations $y_1=0,$  $y_2=0$ and $h^2+4ef=0$). So these two Poisson structures
are not isomorphic, even formally.

The linear Poisson structure corresponding to $\mathfrak{e}(3)$ has the form
$\Pi^{(1)}=x_1\partial x_2\wedge\partial x_3+x_2\partial x_3\wedge\partial
x_1+x_3\partial x_1\wedge\partial x_2+y_1\partial x_2\wedge\partial y_3+ y_2\partial
x_3\wedge\partial y_1+y_3\partial x_1\wedge\partial y_2$ in a natural system of
coordinates. Now put $\Pi=\Pi^{(1)}+\tilde\Pi$ with $\tilde\Pi=(x_1^2+x_2^2+x_3^2)
 (x_1\partial y_2\wedge\partial y_3+x_2\partial y_3\wedge\partial
y_1+x_3\partial y_1\wedge\partial y_2).$
For $\Pi^{(1)}$ the set where the rank is less or equal to 2 is a dimension 3
subspace (given by the equation $y_1 = y_2 = y_3 =0$), while for $\Pi$ the set where
the rank is less or equal to 2 is the origin.

\bibliographystyle{amsalpha}

\providecommand{\bysame}{\leavevmode\hbox to3em{\hrulefill}\thinspace}

\end{document}